\documentclass{SCIS2025}
\numberwithin{equation}{section}
\usepackage{tikz}
\graphicspath{{images/}}
\usepackage{amsmath,amssymb,amsfonts}
\usepackage{graphicx}
\usepackage{hyperref}
\usepackage{textcomp}
\usepackage{lipsum}
\usepackage{amsfonts}
\usepackage{epstopdf}
\usepackage{amsopn}
\usepackage{cite}
\usepackage{amsmath,amssymb,amsfonts}
\usepackage{amsmath,bm}
\usepackage{subfloat}

\begin{document}


\ArticleType{RESEARCH PAPER}
\Year{2025}
\Month{}
\Vol{}
\No{}
\DOI{}
\ArtNo{}
\ReceiveDate{}
\ReviseDate{}
\AcceptDate{}
\OnlineDate{}
\AuthorMark{}
\AuthorCitation{}

\title{On Control Networks Over Finite Lattices}{Ji Z, Cheng D. On control networks over finite lattices}

\author[1,2]{Zhengping JI}{{jizhengping@amss.ac.cn}}
\author[1]{Daizhan CHENG}{}

\address[1]{Key Laboratory of Systems and Control, Academy of Mathematics and Systems Science, \\Chinese Academy of Sciences,
	Beijing 100190, China}
\address[2]{School of Mathematical Sciences, University of Chinese Academy of Sciences, Beijing 100049, China}

\abstract{The modelling and control of networks over finite lattices are studied via the algebraic state space approach. Using the semi-tensor product of matrices, we obtain the algebraic state space representation of the dynamics of (control) networks over finite lattices. Basic properties concerning networks over sublattices and product lattices are investigated, which shows the application of the analysis of lattice structure in the model reduction and control design of networks. Then algorithms are developed to recover the lattice structure from the structural matrix of a network over a lattice, and to construct comparability graphs over a finite set to verify whether a multiple-valued logical network is defined over a lattice. Finally, numerical examples are presented to illustrate the results.}


\maketitle

\section{Introduction}

The Boolean networks were firstly introduced by Kauffman \cite{kau69} to
model the genetic regulatory networks, and were later generalized as multiple-valued logical networks to improve the accuracy of the logical variables. The theory of Boolean algebra proposed by \cite{boo47}, where the value of a logical node is quantized to zero or one, together with the theory of $n$-valued logic introduced by \cite{pos21}, can be unified into general classes of algebras. These notions were generalized in \cite{ros42}  to the so-called Post algebra, which, by revealing itself as a special class of finite lattices, provided a way of treating multiple-valued algebra from the perspective of lattice theory \cite{ep60,tr77}. 

Lattices are common objects in combinatorics, and networks over lattices can be commonly found in modelling and control of communication networks \cite{shu,tu} as well as systems biology \cite{bou}. Since multiple-valued logical networks can be considered as defined over Post algebras which are special cases of lattices, one can in general consider a network where each node takes value in a finite lattice, and the evolution depends on the $\sup$ and $\inf$ operations between the states. These networks are common in applications. For example, distant data verification and breakdown restoration in multi-agent systems with logically linked list of entries and the distributed ledger can be modelled using Allen-Givone algebra, which can be viewed as a special class of lattices \cite{by20,by22}; the fuzzy bisimulation of nondeterministic transition systems is modelled through residuated lattices \cite{qia,qiao}; {the asynchronous dynamics of a gene network with multiple expression levels can be modelled through operators defined on a finite lattice \cite{ric07}. Moreover, if the nodes of networks are defined over lattices such as the finite sets of Eisenstein integers \cite{sun}, revealing the underlying partial order or dimensions of the lattice will benefit the analysis of the network structure. }

On the other hand, networks over lattices have simple logical expressions due to the properties such as commutativity and idempotency of the generating operators, making it mathematically easier to present. From the underlying algebraic structure of the network, one would expect that the following properties hold: restricting such a network to a sublattice provides an invariant sub-network; the control properties of networks over product lattices are completely determined by their decomposed subnetworks over sub-lattices, etc. These properties reduce the complexity of analysis and control design of such systems, as we will show in this paper.

{However, to the authors'  knowledge, there have been few investigations on general (control) networks over lattices. There are two problems as possible reasons for this. First, given the input-output relation of a finite-valued network, we are not always able to recognize whether it is defined over a lattice; secondly, when a network is indeed expressed through the composition of basic operators on a lattice, the expression of the nonlinear dynamics is not easy to be simplified for deriving a transition law, impeding the analysis of the network.}

Consider the first problem. Since networks over lattices may encode information of the underlying partial order structure in its algebraic expression, it is possible that from the dynamics of a given multiple-valued logical network one may recognize the basic operators serving as ``bricks" or ``building blocks" that generates the system, and hence determine if it is actually defined over a lattice. The problem of verifying whether a network can be considered as one defined over a lattice (and further reconstruct this lattice structure if possible) is of practical importance since it is a special case of the NP-hard constraint satisfaction problems \cite{rus10}. Expressing a $k$-valued network over a lattice is not only practically useful but also theoretically challenging. It is well known that even expressing a $3$-valued logical networks into $3$-valued logical form is difficult in general, while the $k$-valued case is even more complicated \cite{ada03}. There have been investigations on the generating systems of function algebras (see, for example, Part II of \cite{lau06}), but as far as the authors are concerned, few research has been done for designing a lattice structure on a finite set to make a (control) network be generated by the basic operators on the lattice, which is an important issue, since sometimes we only have the desired input-output relation of a system and we would like to realize it through simple generating functions \cite{rin77}. For example, if the desired transition law of a network can be written explicitly in terms of operators on a lattice, then the system design can be executed based on formal calculations over these operators. This can be viewed as an discrete-state analogue of the circuit realization problem \cite{gar13,sun06}.

Concerning the second problem, we introduce the semi-tensor product (STP) of matrices as the tool. Since 2009, it has been applied to the study of Boolean networks \cite{che09}, providing a convenient way to represent the dynamics into algebraic equations and promoting the development of the theory of finite-valued networks. Investigations have emerged mainly on two directions. One is the control problems of Boolean networks, such as controllability \cite{che10a}, observability \cite{for13,zhang}, stabilization \cite{li17}, tracking \cite{zha20}, decoupling \cite{liy}, etc. Another direction is extending the values of states in the networks   from the Boolean case (where they are binary) to $k$-valued, mixed \cite{lu17}, finite ring \cite{che21,zou} and finite field cases \cite{men, lin, linl}. {The STP makes it possible to simplify the expression of a finite-valued control network to a unified bilinear from, simplifying the analysis of its control properties. This method has been applied to finite-valued evolutionary games as well \cite{yan}. However, a main obstacle of applying the STP is that the computational complexity increases exponentially with the number of the nodes in the networks; there have been several approaches proposed for this problem, such as aggregation \cite{ji} and pinning control method \cite{zho22}.}

This paper aims to provide a framework for (control) networks over finite lattices to solve the above problems. Using the STP, the algebraic state space representation (ASSR) of such networks is derived, allowing us to analyze the system with existing results in multiple-valued logical networks. We show that when the lattice is the product of some finite lattices, the control  properties are determined by the  subnetworks defined over factor lattices. Further, from a network over a lattice, we recover the underlying lattice structure from its ASSR; when the existence of underlying partial order relations is unknown for an arbitrarily given network, we give necessary conditions for it to be a network over a lattice, and try to construct such a partial order to make the network be generated by some classes of basic operators over a lattice. {Compared with the existing methods of analyzing networks over finite rings \cite{che21} and finite fields \cite{men,lin}, our algorithm for the first time provides a method for reconstructing the algebraic structure over dynamics, while in the biographies, there have been no results derived for recovering the ring or field operations such as addition and multiplication that builds the network dynamics. On the other hand, our algorithm provides a way to reduce the computational complexity of the systems by decomposing the networks into subnetworks over factor lattices. Compared with existing methods for model reduction of logical networks such as aggregation \cite{zha}, it has no restriction on the topological structures of the network, and is applicable for general networks over product lattices.}

The rest of this paper is organized as follows. Section \ref{S2} provides preliminaries about the STP and lattice theory; Section \ref{S3} studies the networks over finite lattices under the framework of algebraic state space and vector expressions: first it gives criteria for a finite set endowed with an operator to be a lattice, and then investigates the basic properties of networks over product lattices; Section \ref{S4} is devoted to constructing and recovering the underlying lattice structures for finite-valued networks, designing algorithms to construct partial order so that a given network can be viewed as generated by basic operators of a lattice. These results are illustrated by numerical examples in Section \ref{S5}. Section \ref{S6} consists of conclusions and further problems.

Before ending this introduction we give a list for notations which are used in the sequel:

\begin{enumerate}

 	\item ${\mathcal M}_{m\times n}$: set of $m\times n$-dimensional real matrices.
	
	\item $\mathrm{Col}(A)$ ($\mathrm{Row}(A)$): the set of columns (rows) of $A$; $\mathrm{Col}_i(A)$ ($\mathrm{Row}_i(A)$): the $i$-th column (row) of $A$.
	
	\item $\delta_n^i$: $i$-th column of the identity matrix $I_n$.
	
	\item ${\mathcal D}_k:=\{0,\cdots,k-1\}$.
	
	\item $\Delta_k:=\mathrm{Col}(I_k)$.

	\item ${\mathcal L}_{m\times n}$: the set of logical matrices (A matrix $L\in {\mathcal M}_{m\times n}$ is called a logical matrix, if $\mathrm{Col}(L)\subset \Delta_m$).
	
	\item $\delta_m[i_1,\cdots,i_n]$: a brief notation for logical matrices, that is,  $\delta_m[i_1,\cdots,i_n]:=\left[\delta_m^{i_1},\cdots,\delta_m^{i_n}\right]$.
	
	\item $A \times_{\mathcal B} B$: the Boolean product of $A\in B_{m\times n} $, $B \in B_{n\times p}$, i.e., $[A \times_{\mathcal B}B ]_{i,j} = 0$ if $[AB]_{i,j} = 0$, and $[A \times_{\mathcal B}B ]_{i,j} = 1$ if $[AB]_{i,j} > 0$, $i=1,\cdots,m$, $j=1,\cdots,p$.
	
\end{enumerate}

\section{Preliminaries}\label{S2}

\subsection{STP and Algebraic Expression of Multiple-Valued Networks}
We first give a brief review on STP of matrices which is the main tool in this paper. We refer to \cite{che12} for details.

\begin{definition}\label{d2.1.1} Let $A\in {\mathcal M}_{m\times n}$,  $B\in {\mathcal M}_{p\times q}$ be real matrices, and the least common multiple of $n$ and $p$ be $t=\mathrm{lcm}\{n,p\}$.
	The STP of $A$ and $B$, denoted by $A\ltimes B$, is defined as
	\begin{equation} \label{2.1.1} \left(A\otimes I_{t/n}\right)\left(B\otimes I_{t/p}\right).
	\end{equation}
	where $I_k$ is the $k\times k$ identity matrix, and $\otimes$ is the Kronecker product.
\end{definition}

Throughout this paper, all products are assumed to be semi-tensor products and the symbol $\ltimes$ is usually omitted.

When the dimensions of two matrices are compatible, the STP is the same as the conventional matrix product, and the laws concerning associativity, distributivity, transpose and inverse hold for it as well. Further, it has the following properties concerning commutativity.
\begin{proposition}[\cite{che12}]\label{p2.1.5}
	Let $X\in \mathbb{R}^m$ be a column vector and $M$ be a matrix. Then $X\ltimes M=\left(I_m\otimes M\right)\ltimes X$.
 
	Given two column vectors ~$X\in \mathbb{R}^m$, $Y\in \mathbb{R}^n$, then
	\begin{align} \label{2.1.6} 
		W_{[m, n]}\ltimes X\ltimes Y=Y\ltimes X
	\end{align}
	where $W_{[m,n]}\in {\mathcal M}_{mn\times mn}$ is the $(m,n)$-order swap matrix defined as
	\begin{align*}
			W_{[m,n]}=\delta_{mn}[1,m+1,\cdots,(n-1)m+1,2,m+2,\cdots,(n-1)m+2,\cdots,m,2m,\cdots,nm].
	\end{align*}
\end{proposition}

\begin{definition}\label{d2.1.6}
	\begin{enumerate}
		\item Let $x_i\in {\mathcal D}_{k_i}$, $i=1,\cdots,n$. A map $f: \prod_{i=1}^n{\mathcal D}_{k_i}\rightarrow {\mathcal D}_{k_0}$ is called a multiple-valued logical function.
		If $k_1=k_2=\cdots=k_n=k_0$, $f$ is called a $k_0$-valued logical function.
		\item Let $x_i(t)\in {\mathcal D}_{k_i}$, $f_i:\prod_{i=1}^n{\mathcal D}_{k_i}\rightarrow {\mathcal D}_{k_i}$, $i=1,\cdots,n$. The system
		\begin{align}\label{2.1.8}
			\begin{cases}
				x_1(t+1)=f_1(x_1(t),\cdots,x_n(t))\\
				\vdots~~\\
				x_n(t+1)=f_n(x_1(t),\cdots,x_n(t)),\\
			\end{cases}
		\end{align}
		is called a multiple-valued logical dynamic system.

		If $k_1=k_2=\cdots=k_n=k_0$, the system is called a $k_0$-valued logical dynamic system.
	\end{enumerate}
\end{definition}

Identify $i\sim \delta_k^{i}, ~i=1,\cdots,k-1$, and $0\sim\delta_k^k$, then $x\in {\mathcal D}_k$ can be expressed as $x\in \Delta_k$. The latter is called the vector expression of a logical variable. Using this and the above properties of STP, one can express a multiple-valued logical function and further a network in the algebraic state space.

\begin{theorem}[\cite{che11}]\label{t2.1.7}
	Given a multiple-valued logical function $f: \prod_{i=1}^n{\mathcal D}_{k_i}\rightarrow {\mathcal D}_{k_0}$ denoted as $y=f(x_1,\cdots,x_n)$, there exists a unique $M_f\in {\mathcal L}_{k_0\times k}$, where $k=\prod_{i=1}^nk_i$, such that in vector expression we have
	\begin{align*}
		y=M_f\ltimes_{i=1}^n x_i.
	\end{align*}
	$M_f$ is called the structure matrix of $f$.
\end{theorem}

Applying Theorem \ref{t2.1.7} to each equation of (\ref{2.1.8}), we have
\begin{align}\label{2.1.11}
	\begin{cases}
		x_1(t+1)=M_1\ltimes_{i=1}^nx_i(t)\\
		\vdots~~\\
		x_n(t+1)=M_n\ltimes_{i=1}^nx_i(t),\\
	\end{cases}
\end{align}
where $M_i$ is the structure matrix of $f_i$, $i=1,\cdots,n$.

\begin{theorem}[\cite{che11}]\label{t2.1.10} Denote by $x(t)=\ltimes_{i=1}^nx_i(t)$. Then (\ref{2.1.11}) can be expressed as
	\begin{align}\label{2.1.13}
		x(t+1)=Mx(t),
	\end{align}
	where $M=M_1*\cdots*M_n$, and $*$ is the Khatri-Rao product of matrices. (\ref{2.1.13}) is called the algebraic state space representation (ASSR) of the multiple-valued logical network (\ref{2.1.8}), and $M$ is called its structure matrix.
\end{theorem}

\subsection{Finite Lattices}
Next, we review the basic notions of lattice theory.

A lattice is a partially ordered set $(L,\leqslant)$ where each pair of elements $x,y\in L$ has a least upper bound (denoted by $x\vee y$ or $\sup(x,y)$) and a greatest lower bound (denoted by $x\wedge  y$ or $\inf(x,y)$). {For example, define a partial order on $\mathbb{N}$ as $a\leqslant b\Leftrightarrow a|b$, and let $a\vee b$ be the least common multiple and $a\wedge  b$ the greatest common divisor of two integers $a$, $b$, then $(\mathbb{N},\leqslant)$ is a lattice. Let $L:=\{0,1\}$, then $L$ is a lattice under natural order if we define the operators $\vee$, $\wedge $ as disjunction and conjunction of the Boolean variables respectively.}

An alternative definition is that a lattice is a tuple $(L,\wedge,\vee)$ where $\wedge,\vee$ are binary operators on $L$ satisfying commutativity, associativity, idempotency and the absorption property. Since only finite lattices are considered in this paper, we adopt the following definition.

\begin{definition}[\cite{ber15}]\label{d3.1}
	A finite set $L$ is called a lattice, if there exists a binary operator $\vee$ on $L$, satisfying
	\begin{enumerate}
		\item[(i)] $x\vee x=x$;
		\item[(ii)] $x\vee y=y\vee x$;
		\item[(iii)] $(x\vee y)\vee z=x\vee (y\vee z)$;
		\item[(iv)] $\exists \mathbf{0}\in L$, s.t. $w\vee\mathbf{0}=w$, $\forall w\in L$,
	\end{enumerate}
	where $x, y, z$ are arbitrary elements in $L$.
\end{definition}

\begin{remark}
	If we define a new binary operator $\wedge$ on $L$ as $a\wedge b:=\bigvee_{u\in S}u$, where $S:=\{u\in L|u\leqslant a,u\leqslant b\}$, then the triple $(L,\wedge,\vee)$ coincides with the conventional definition of a lattice, which means there exists a partial order $\leqslant$ on $L$ such that $\vee$, $\wedge$ are least upper bound and greatest lower bound operators under this order, and $(L,\leqslant)$ admits unique maximal and minimal elements. The proofs can be found in \cite{ber15,sta12}.
\end{remark}

\begin{figure}
	\centering
	\setlength{\unitlength}{0.3 cm}
	\begin{picture}(6,6)
		\thicklines
		\put(3,1){\circle*{0.2}}
		\put(1,3){\circle*{0.2}}
		\put(5,3){\circle*{0.2}}
		\put(3,5){\circle*{0.2}}
		\put(1,3){\line(1,-1){2}}
		\put(1,3){\line(1,1){2}}
		\put(5,3){\line(-1,-1){2}}
		\put(5,3){\line(-1,1){2}}
		\put(2.9,5.2){\small$P_1$}
		\put(2.9,0.2){\small$P_4$}
		\put(-0.2,3.3){\small$P_2$}
		\put(5.2,3){\small$P_3$}
	\end{picture}
	\caption{Hasse diagram of a 4-element lattice $L$\label{Fig.3.1}}
\end{figure}
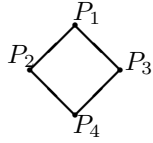

{A lattice $L$, like other partially ordered sets, can be completely characterized by its Hasse diagram \cite{sta12}: a graph whose vertices are the elements in the lattice. In the Hasse diagram of a lattice, there is an edge between two vertices $a,b$ if and only if the elements represented by these two vertices are comparable ($a\leqslant b$ or $b\leqslant a$), and if $b$ covers $a$, which means $a\leqslant b$ and there is no element $c\in L$ such that $a\leqslant c\leqslant b$, then $b$ is drawn ``above" $a$ (i.e., with a higher vertical coordinate).}

{For example, consider a lattice $L=\{P_1,P_2,P_3,P_4\}$ whose Hasse diagram is shown as in Figure \ref{Fig.3.1}. One can see that $P_1$ is the greatest element and $P_4$ the least, while $P_2$ and $P_3$ are incomparable, $P_2\vee P_3=P_1$, $P_2\wedge P_3=P_4$.
}

\section{Networks Over Finite Lattices}\label{S3}

We first use the STP to give an example of how algebraic state space will help to solve the problems in lattice theory.

Since an operator over a finite set is a logical function (according to Definition \ref{d2.1.6}), we can easily adapt the STP framework to networks over finite lattices. For example, the following proposition gives criteria for a finite set endowed with a binary operator to be a lattice.
\begin{proposition}\label{p3.3}
	Let $L$ be a finite set of cardinal $k$ and $f$ a binary operator over $L$ with a structure matrix $M=[M_1,\cdots,M_k]$ where $M_i\in L_{k\times k}$, $i=1,\cdots, k$. Then $(L,f)$ is a lattice if and only if $M$ satisfies the following conditions:
	\begin{enumerate}
		\item[(i)] $M\mathrm{diag}(\delta_k^1,\cdots,\delta_k^k)=I_k$;
		\item[(ii)] $MW_{[k,k]}=M$;
		\item[(iii)] $M^2=M(I_k\otimes M)$;
		\item[(iv)] $\exists i\in {\mathcal D}_k$, s.t. $M_i=I_k$. 
	\end{enumerate}
\end{proposition}
The proof  follows directly from the basic properties of STP stated in the previous section. This gives a convenient way to verify through the structure matrix whether a binary function over a finite set defines a lattice.

\subsection{ASSR of Control Networks Over Lattices}

Let $F:L^n\rightarrow L$ be an $n$-ary function on $L$. If $F$ is obtained from the {composition} of $\vee$ and $\wedge$ on $L$, we call it a lattice function. {That is to say, up to inserting parentheses, $F$ is written in the following form 
$$
F(x_1,\cdots,x_n)=x_{i_1}\bigcirc x_{i_2}\bigcirc\cdots\bigcirc x_{i_N},\quad i_1,\cdots,i_N\in\{1,\cdots,n\}, \quad N>0,
$$
where the operator $\bigcirc$ is either $\vee$ or $\wedge$.}
Since lattice functions are special cases of multiple-valued logical functions, we may apply the STP method to model the networks over finite lattices.

In general, a network over a finite lattice $L$ has the form of (\ref{2.1.8}) where $x_i(t)\in L$ are state variables and the maps $f_i:L^n\rightarrow L$ are lattice functions, $i=1,\cdots,n$. When $L={\mathcal D}_2$, the system becomes a standard Boolean network.

When there are controls $u_1(t), \cdots,u_m(t)\in L$ in the network, its dynamics can be expressed as
\begin{small}
	\begin{align}\label{ev.1.4}
		\begin{cases}
			x_1(t+1)=g_1(x_1(t),\cdots,x_n(t);u_1(t),\cdots,u_m(t))\\
			\vdots\\
			x_n(t+1)=g_n(x_1(t),\cdots,x_n(t);u_1(t),\cdots,u_m(t)),\\
		\end{cases}
	\end{align}
\end{small}
where $g_i$ are lattice functions, $i=1,\cdots,n$.

Assume that $\mathrm{card}(L)=\kappa$. Since the systems over finite lattices belong to the class of multiple-valued logical systems, by Theorem \ref{t2.1.10}, using STP, (\ref{2.1.8}) can be converted into its algebraic state space representation as
\begin{align}\label{ev.1.5}
	x(t+1)=Mx(t),
\end{align}
where $x(t)=\ltimes_{j=1}^nx_j(t)$, $M\in {\mathcal L}_{\kappa^n \times \kappa^n}$, and (\ref{ev.1.4}) as
\begin{align}\label{ev.1.6}
	x(t+1)=Pu(t)x(t),
\end{align}
where $u(t)=\ltimes_{s=1}^mu_s(t)$, $P\in {\mathcal L}_{\kappa^n \times \kappa^{n+m}}$.

The STP provides a way to rewrite the nonlinear control system (\ref{ev.1.4}) to a simplified and unified form, giving it a mathematically neat expression. We give a pedagogical example to illustrate how ASSR is derived.

\begin{example}\label{5.0}
	Consider a $4$-element lattice $L$ with Hasse diagram as shown in Figure \ref{Fig.3.1}. Setting $P_i\sim \delta_4^i$, $i=1,2,3,4$, {and using the formula of the operators $\vee$, $\wedge$ derived from the diagram as
    $$
    \begin{array}{l}
         P_1\vee P_2=P_1, ~P_1\vee P_3=P_1, ~P_1\vee P_4=P_1, ~P_2\vee P_3=P_1, ~P_2\vee P_4=P_2, ~P_3\vee P_4=P_3, \\
         P_1\wedge P_2=P_2, ~P_1\wedge P_3=P_3, ~P_1\wedge P_4=P_4, ~P_2\wedge P_3=P_4, ~P_2\wedge P_4= P_4, ~P_3\wedge P_4=P_4,
    \end{array}
    $$}
    we can solve the structure matrix of the operators $\vee$, $\wedge$ respectively {(for details of calculating the structure matrices of finite-valued functions, one may refer to \cite{che11})}, that is, $\forall x,y\in\{P_1,\cdots,P_4\}$,
	\begin{align*}
		&x\vee y=M_{\vee}xy, \quad x\wedge y=M_{\wedge}xy,\\
		&M_{\vee}=\delta_4[1,1,1,1,1,2,1,2,1,1,3,3,1,2,3,4];\\
		&M_{\wedge}=\delta_4[1,2,3,4,2,2,4,4,3,4,3,4,4,4,4,4].
	\end{align*}
	Here we make no distinction between an element in the lattice and its vector expression. Next, assume a network over $L$ is defined as
	\begin{align}\label{1}
		\begin{cases}
			x_1(t+1)=x_1(t)\vee (x_2(t)\wedge u(t)),\\
			x_2(t+1)=x_1(t)\vee x_2(t)\vee u(t).
		\end{cases}
	\end{align}
	
	Then using the commutative properties in Proposition \ref{p2.1.5}, it is easy to construct the following component-wise ASSR from the structure matrices of the operators.
	$$
	\begin{cases}
		x_1(t+1)=M_1u(t)x(t),\\
		x_2(t+1)=M_2u(t)x(t),
	\end{cases}
	$$
	where $x(t)=x_1(t)x_2(t)$, and
	\begin{small}
		$$
		\begin{array}{l}
			M_1=M_{\wedge}(I_4\otimes M_{\vee})W_{[4,16]},\\
			M_2=M_{\vee}(I_4\otimes M_{\wedge})W_{[4,16]}.
		\end{array}
		$$
	\end{small}
	Then by a construction similar with (\ref{2.1.13}), the ASSR is $x(t+1)=Mu(t)x(t)$, where
	\begin{small}
		$$
		\begin{array}{ccl}
			M&=&M_1*M_2\\
			~&=&\delta_{16}[1,1,1,1,5,6,5,6,9,9,11,11,13,14,15,16,1,5,1,5,6,6,6,6,9,13,11,15,14,14,16,16,1,1,\\
			~&~&9,9,5,6,13,14,11,11,11,11,15,16,15,16,1,5,9,13,6,6,14,14,11,15,11,15,16,16,16,16].
		\end{array}
		$$
	\end{small}

\end{example}
\begin{remark}
	Note that one advantage of networks over lattices is that they are composed of basic binary operators so that we may use STP to solve the structure matrix solely from its operator expression, which cannot be done for arbitrary multiple-valued networks. The ASSR make the computation and analysis of system (\ref{ev.1.4}) much easier, since it allows one to use existing results in multiple-valued logical networks to investigate the control properties of the system. For example, one may consider problems such as controllability and observability of (\ref{1}).
\end{remark}

\begin{example}
	Consider the controllability problem of the network (\ref{1}). Denote the matrix $M$ by $[M_1,\cdots,M_{16}]$ where $M_i\in {\mathcal L}_{16\times16}$, following \cite{che11}, construct its controllability matrix as
	\begin{align}\label{3.1}
		C=\sum_{j=1}^{16}\big(\sum_{i=1}^{16}M_i\big)^{(j)},
	\end{align}
	Taking the addition and multiplication in (\ref{3.1}) to be Boolean type (that is to say, $1+1=1$), we solve the controllability matrix as
	\begin{align*}	
		\begin{tiny}
			C=\left(\begin{array}{cccccccccccccccc}
				1 & 1 & 1 & 1 & 0 & 0 & 0 & 0 & 0 & 0 & 0 & 0 & 0 & 0 & 0 & 0 \\
				0 & 0 & 0 & 0 & 0 & 0 & 0 & 0 & 0 & 0 & 0 & 0 & 0 & 0 & 0 & 0 \\
				0 & 0 & 0 & 0 & 0 & 0 & 0 & 0 & 0 & 0 & 0 & 0 & 0 & 0 & 0 & 0 \\
				0 & 0 & 0 & 0 & 0 & 0 & 0 & 0 & 0 & 0 & 0 & 0 & 0 & 0 & 0 & 0 \\
				0 & 1 & 0 & 1 & 1 & 0 & 1 & 0 & 0 & 0 & 0 & 0 & 0 & 0 & 0 & 0 \\
				0 & 1 & 0 & 1 & 1 & 1 & 1 & 1 & 0 & 0 & 0 & 0 & 0 & 0 & 0 & 0 \\
				0 & 0 & 0 & 0 & 0 & 0 & 0 & 0 & 0 & 0 & 0 & 0 & 0 & 0 & 0 & 0 \\
				0 & 0 & 0 & 0 & 0 & 0 & 0 & 0 & 0 & 0 & 0 & 0 & 0 & 0 & 0 & 0 \\
				0 & 0 & 1 & 1 & 0 & 0 & 0 & 0 & 1 & 1 & 0 & 0 & 0 & 0 & 0 & 0 \\
				0 & 0 & 0 & 0 & 0 & 0 & 0 & 0 & 0 & 0 & 0 & 0 & 0 & 0 & 0 & 0 \\
				0 & 0 & 1 & 1 & 0 & 0 & 0 & 0 & 1 & 1 & 1 & 1 & 0 & 0 & 0 & 0 \\
				0 & 0 & 0 & 0 & 0 & 0 & 0 & 0 & 0 & 0 & 0 & 0 & 0 & 0 & 0 & 0 \\
				0 & 0 & 0 & 1 & 0 & 0 & 1 & 0 & 0 & 1 & 0 & 0 & 1 & 0 & 0 & 0 \\
				0 & 0 & 0 & 1 & 0 & 0 & 1 & 1 & 0 & 1 & 0 & 0 & 1 & 1 & 0 & 0 \\
				0 & 0 & 0 & 1 & 0 & 0 & 1 & 0 & 0 & 1 & 0 & 1 & 1 & 0 & 1 & 0 \\
				0 & 0 & 0 & 1 & 0 & 0 & 1 & 1 & 0 & 1 & 0 & 1 & 1 & 1 & 1 & 1 
			\end{array}\right),
		\end{tiny}
	\end{align*}
	which gives full information of controllability of each state in (\ref{1}). That is to say, a state $\delta_{16}^i$ is controllable from $\delta_{16}^j$ if and only if the $(i,j)$-th entry of $C$ is nonzero (for details, one may refer to \cite{zh}).
	
\end{example}

Another special property for networks over lattices is that the subnetworks defined over a sublattice of $L$ are invariant subspaces of the network over $L$, that is, the dynamics starting from a point in the sublattice will remain in it. 

\begin{example}
	Recall the system (\ref{1}) and remove the controls in the network. Assume $S=\{P_1,P_2,P_4\}\subset L$, which is clearly a sub-lattice of $L$.
	Re-assign the nodes as $P_1\sim \delta_3^1$, $P_2\sim\delta_3^2$, $P_4\sim\delta_3^3$. It is easy to figure out that setting $z=x|_S$  we have
	$$
	z(t+1)=M|_S z(t),
	$$
	where $M|_S=\delta_9[1,1,1,4,5,5,7,8,9]$.
	
	If we add the controls back and restrict them to the sublattice $S$, then it turns out that the system is a control-invariant subspace which can be viewed as defined over the sublattice with the dynamics 
	$$
	z(t+1)=P|_Su(t)z(t),
	$$
	where
	$$
		P|_S=\delta_9[1,1,1,2,5,5,3,6,9,1,2,2,5,5,5,6,6,9,1,2,3,5,5,6,9,9,9].
	$$
\end{example}

\subsection{Control Networks Over Product Lattices}\label{S3B}

Next we consider the networks over a special class of finite lattices, which are the Cartesian product of some finite lattices. Suppose $L=L_1\times\cdots\times L_p$, where $L_i$ is a finite lattice of cardinal $k_i$, $i=1,\cdots,p$. Define the partial order relation $\leqslant$ on $L$ by
$$
(x_1,\cdots,x_p)\leqslant(y_1,\cdots,y_p)\Leftrightarrow x_1\leqslant y_1,\cdots,x_p\leqslant y_p.
$$
Then a network over $L$ with ASSR (\ref{ev.1.6}), denoted by $\Sigma$, can be decomposed into several factors $\Sigma_1, \cdots,\Sigma_p$, where the factor $\Sigma_i$ is a network consisting of the corresponding $i$-th factor of the state $(x_1(t),\cdots,x_p(t))$, with the binary operations in $L_i$, $i=1,\cdots,p$.

Then a decomposition theorem follows.

\begin{theorem}\label{t3}
	The system $\Sigma$ is controllable (observable, synchronizable, stabilizable), if and only if, the factor systems $\Sigma_i$, $i=1,\cdots,p$ are controllable (observable, synchronizable, stabilizable).
\end{theorem}
\noindent {\it Proof.}
We only need to show that a network defined over a product lattice of $L_1$ and $L_2$ can be decomposed into two systems $\Sigma_1$, $\Sigma_2$, which are defined over $L_1$, $L_2$ respectively, i.e. any trajectory $\{z(t)|t\geqslant0\}\subset L$ starting from a point $(x_0,y_0)$ can be decomposed into two trajectories $\{x(t,x_0)|t\geqslant0\}$, $\{y(t,y_0)|t\geqslant0\}$ starting from points $x_0\in L_1$, $y_0\in L_2$ respectively, such that $z(t)=(x(t,x_0),y(t,y_0))$, $\forall t>0$. Since a factor of a product lattice can be viewed as a sublattice, the factor lattices are both invariant subspaces, such decomposition is obvious, and the conclusion follows.
\hfill $\Box$
{
\begin{remark}
	For a network over the above product lattice with an ASSR (\ref{ev.1.5}), the ASSR of its subsystem over the factor lattice $L_i$, $i=1,\cdots,p$ is derived as
	$$
	\tilde{x}^i(t+1)=H_i\times_{\mathcal B}M\times_{\mathcal B} H_i^{\mathrm{T}}\tilde{x}^i(t),
	$$
	where $\tilde{x}^i(t)\in \Delta_{k_i^n}$, and  $H_i\in {\mathcal L}_{k_i^n\times k^n}$ is the structure matrix of the projection map $L\rightarrow L_i$. This is due to the fact that $L_i$ is invariant with respect to the evolution of the subnetwork $\Sigma_i$. For details of the construction of the above structure matrix, one may refer to \cite{ji}.
\end{remark}
}

{ 
    As is known, when we adopt the STP method, the complexity of analyzing the logical network increases exponentially with respect to the number of nodes \cite{che09,che10a}. Applying Theorem \ref{t3}, we may reduce the complexity of analyzing the network over a product lattice $L=L_1\times\cdots\times L_p$ by decomposing it into systems over subsystems over sublattices $L_i$,  $i=1,\cdots,p$  and reducing its dimension. For example, suppose that the network has $n$ nodes, $\mathrm{card}(L_i)=k_i$, and $k=k_1\cdots k_p$. Then the dimension of the ASSR of the system is $k^n$. If the complexity of executing an algorithm for analysis of the network is $O(k^{n})$, then after the decomposition, the complexity will be reduced to $O(\sum_{i=1}^pk_i^n)$, which greatly eases the computation load.
}

{For an application of Theorem \ref{t3}, we give the an example of systems evolving over product lattices. Consider the following linear switched system with logical switching \cite{gu}:}
\begin{subequations}\label{e:parent}
\begin{align}
    &x(t+1)=A_{\sigma(t)}x(t)+B_{\sigma(t)}u(t),\label{e:1}\\
    &\eta(t)=Q(x(t))\label{e:2},\\
     &  \begin{cases}
        \lambda(t+1)=\gamma(\lambda(t),\eta(t)),\\
        \sigma(t)=\beta(\lambda(t),\eta(t)),
    \end{cases}\label{e:3}
\end{align}
    \end{subequations}
where $Q:D\rightarrow {\mathcal D}_k$ is a quantizer, $D=[\eta_1,\xi_1)\times\cdots[\eta_n,\xi_n)\subset \mathbb{R}^n$, $\sigma(t)\in{\mathcal D}_N$ is the switching signal, (\ref{e:3}) is an $\ell$-valued logical dynamic system with $\gamma$, $\beta$ logical functions and $\lambda(t)\in {\mathcal D}_{\ell}$ a logical variable, $k$, $N$, $\ell$ are integers,  $A_i$, $B_i$ are matrices, $\eta_i<\xi_i$ are real numbers, $i=1,\cdots,N$. $x(t)\in D$ is the continuous-state variable. It is common that the the sampling of the switching signal is based on coordinate partitions \cite{pol}, that is to say, the quantizer $Q$ is defined as
\begin{align*}
    Q(x_1,\cdots,x_n)=\alpha_{i_1,i_2,\cdots,i_n}, ~\text{if}~\beta^1_{i_1}\leqslant x_1 <\beta^1_{i_1+1},~\beta^2_{i_2}\leqslant x_2<\beta^2_{i_2+1},\cdots,~\beta^n_{i_n}\leqslant x_n<\beta^n_{i_n+1},
\end{align*}
where $\alpha_{i_1,\cdots,i_n}\in{\mathcal D_k}$ pairwise distinct over $i_1=1,\cdots,N_1$, $i_2=1,\cdots,N_2$, $\cdots$, $i_n=1,\cdots,N_n$, with $N_1,\cdots,N_n$ integers. In short, the quantizer assigns each state $x=(x_1,\cdots,x_n)\in\mathbb{R}^n$ an integer in ${\mathcal D}_k$, according to the interval $[\beta^j_{i_j},\beta^j_{i_j+1})$ that each component $x_j$ lies in, and the domain $D\subset\mathbb{R}^n$ is partitioned as
$$
D=\Big\{[\beta^1_1,\beta^1_2)\cup[\beta^1_2,\beta^1_3)\cup\cdots\cup[\beta^1_{N_1},\beta^1_{N_1+1})\Big\}\times\cdots\times\Big\{[\beta^n_1,\beta^n_2)\cup[\beta^n_2,\beta^n_3)\cup\cdots\cup[\beta^n_{N_n},\beta^n_{N_n+1})\Big\}.
$$
Apparently, such partition of $D\subset \mathbb{R}^n$ gives rise to a product lattice structure over ${\mathcal D}_k$, and ${\mathcal D}_k$ can be hence decomposed into components $D_{N_1}\times\cdots\times{\mathcal D}_{N_n}$. Since the discretization of the linear switched system (\ref{e:1}) with respect to the quantizer (\ref{e:2}) is a finite transition system over ${\mathcal D}_k$, 
{ if the dynamics of the discretized transition system of $\eta(t)$ is expressed through the operators over this lattice, then the system after discretization becomes a dynamics evolving on a product lattice. The analysis of system (\ref{e:parent}) relies on the mergence of the discretized variable $\eta(t)$ and the logical variable $\lambda(t)$ \cite{gu}, therefore when $\eta(t)$ evolves over a product lattice, using the decomposition of the system over factor sublattices, the analysis for the merged system can be simplified.}

However, the full potential of Theorem \ref{t3} has not yet been revealed if the underlying lattice structure of the network is known a priori, since in this case one can see from the beginning that the system is built up from subsystems over sublattices. In the next section we will focus on (control) networks whose lattice structure is unknown, and Theorem \ref{t3} will be useful for those systems which are found out by our algorithms to be ones defined over product lattices.

\section{Recovering the Lattice Structures}\label{S4}

After presenting the general expressions and control property analysis of  the networks over lattices, we consider the verification of such networks. As aforementioned, networks over lattices are special cases of multiple-valued logical networks; so a natural question rises: how can one recognize them from the vast class of networks over finite sets, and further, how can one construct or recover the order relationships on the underlying sets from the algebraic state space expressions of networks?

The aim of this section is twofold: first, to reconstruct the lattice structure on ${\mathcal D}_k$ from the structural matrix of a network defined over a lattice of $k$ elements; second, to verify if there exists a partial order on ${\mathcal D}_k$ allowing a given network to be one defined over a lattice, or at least be generated by some basic operators on the lattice.

\subsection{Recovering Order Relations Through Lattice Functions}

Consider a system over ${\mathcal D}_k$ with ASSR (\ref{ev.1.5}). The following algorithm can decompose it into the component-wise expression, The proof of which is straightforward.

\begin{proposition}\label{p4.1}
	Consider a matrix $M\in{\mathcal L}_{k^n\times k^n}$, there exists a unique decomposition $M=M_1*\cdots *M_n$, where $M_i=R_iM\in{\mathcal L}_{k\times k^n}$,
	and $R_i={\mathbf 1}_{k^{i-1}}\otimes I_k\otimes{\mathbf 1}_{k^{n-i}}$, ${\mathbf 1}_n$ is the $n$-dimensional row vector with all entries equal $1$.
\end{proposition}

In the following we will show that the information of comparability between the elements in a lattice can be encoded in the structure matrix of any lattice function over it. That is to say, if $M$ is the structure matrix of an $n$-node network defined over a lattice, then any $M_i$ constructed as in Proposition \ref{p4.1}, which is the structure matrix of the dynamics of the $i$-th node, $i=1,\cdots,n$, will be enough to rebuild the lattice structure.

The case of $2$-ary functions, which corresponds to a network of $2$ nodes, is trivial, since by the so-called absorption property a binary lattice function can only be $\vee$ or $\wedge$ (after some procedure of simplification). Therefore it suffices to check the conditions in Proposition \ref{p3.3} on its structure matrix, and recover the partial order by $x\vee y=y\Leftrightarrow x\leqslant y$.

The general $n$-ary case requires a different approach. Without loss of generality, assume that functions discussed in the following do not have ``dumb" indices, i.e.
\begin{align}\label{4.111}
    \begin{array}{l}
        \forall i=1,\cdots,n, ~\exists a_i\neq b_i\in{\mathcal D}_k, ~ \text{s.t.}\\
	f(x_1,\cdots, x_{i-1},a_i,x_{i+1},\cdots,x_n)\neq f(x_1,\cdots, x_{i-1},b_i,x_{i+1},\cdots,x_n). 
    \end{array}
\end{align}
{Specifically, in the network (\ref{2.1.8}), if the index $i$ is dumb in the function $f_j$, it means that the variable $x_i$ does not influence the one-step transition of the variable $x_j$.} In fact, one can verify if a variable in the expression of a multiple-valued logical function is dumb by using the method provided by \cite{che10a}, and further remove them from the explicit expression of the function.

First we point out the fact that given two comparable elements $a$, $b$ in any lattice $(L,\vee,\wedge)$, the operators $\vee$, $\wedge$ restricted to $\{a,b\}$ behave in the same way as the disjunction and conjunction operators over the classical Boolean algebra $\{0,1\}$. Hence we may try to use the ``restricted" structure matrix of a $k$-valued logical function to find all comparable pairs on ${\mathcal D}_k$.

For an $n$-ary logical function $f$ over ${\mathcal D}_k$, consider its restriction to a pair of elements $a,b\in {\mathcal D}_k$, that is to say, the value $f(x_1,\cdots,x_n)$ takes when $x_1,\cdots,x_n$ take values in $\{a,b\}$ only. This transforms $f$ to another logical function $f_{ab}:{\mathcal D}_2^n\rightarrow {\mathcal D}_k$.

Denote by $M\in{\mathcal L}_{k\times k^n}$ the structure matrix of $f$. $\forall a, b\in {\mathcal D}_k$, if their vector forms are $\delta_k^{i_a}$, $\delta_k^{i_b}$, then identify them as $\delta_2^1$, $\delta_2^2$ respectively. By making no distinction between a variable and its vector form ($a\sim\delta_k^{i_a}$, $b\sim\delta_k^{i_b}$, $1\sim\delta_2^1$, $0\sim\delta_2^2$), we may express a variable $x\in\{a,b\}$ through a Boolean one $\tilde{x}\in\{0,1\}$ in vector form as
\begin{align}\label{4.00}
	x=\delta_k[i_a,i_b]\tilde{x}
\end{align}
hence by substitution, the restriction of $f$ on $\{a,b\}^n$, viewed as a function ${f}_{ab}:{\mathcal D}_2^n\rightarrow{\mathcal D}_k$, has vector form expression as
\begin{align}\label{4.01}
	f(x_1,\cdots,x_n)\big|_{\{a,b\}^n}={f}_{ab}(\tilde{x}_1,\cdots,\tilde{x}_n)={M}_{ab}\tilde{x}_1\cdots \tilde{x}_n,
\end{align}
where
\begin{align}\label{4.1}
	{M}_{ab}=M\prod_{j=0}^{n-1}(I_{2^j}\otimes \delta_k[i_a,i_b]),
\end{align}
The structure matrix (\ref{4.1}) is obtained from (\ref{4.00}) and the commutative properties in Proposition \ref{p2.1.5}.

\begin{theorem}\label{t4.2}
	{Let ${\mathcal D}_k$ endowed with a  partial order relation $\leqslant$ be a  lattice, and $f:{\mathcal D}_k^n\rightarrow {\mathcal D}_k$ be a lattice function satisfying (\ref{4.111}).} Given two elements $a,b\in{\mathcal D}_k$ with vector expressions $\delta_k^{i_a}$, $\delta_k^{i_b}$ respectively, $a$ and $b$ are comparable ($a\leqslant b$ or $b\leqslant a$), if and only if,
	\begin{subequations}\label{eq:parent}
		\begin{align}
			&\mathrm{Col}(M_{ab})\subset\{\delta_k^{i_a},\delta_k^{i_b}\};\label{eq:s1}\\
			&\mathrm{Col}(M_{ab}^i)\subset\{\delta_k^{i_a}-\delta_k^{i_b},\mathbf{0}\},~i=1,\cdots,n,\label{eq:s2}\\
			&M_{ab}\delta_{2^n}^1=\delta_k^{i_a},~M_{ab}\delta_{2^n}^{2^n}=\delta_k^{i_b},\label{eq:s3}
		\end{align}
	\end{subequations}
	where $\mathbf{0}=(0,\cdots,0)\in \mathbb{R}^k$, $M_{ab}$ is defined as in (\ref{4.1}), and
	\begin{align*}
		M_{ab}^i:=M_{ab}W_{[2,2^{i-1}]}\begin{pmatrix}
			1\\-1
		\end{pmatrix}.
	\end{align*}
\end{theorem}
\noindent {\it Proof.}
(Necessity) When $a$, $b$ are comparable, since $\vee$, $\wedge$ restricted to $\{a,b\}$ are identical to disjunction and conjunction operators on the classical Boolean algebra $\{0,1\}$ (denoted by $\vee_2$, $\wedge_2$ respectively), $f|_{\{a,b\}^n}$ as a lattice function can be viewed as belonging to the class generated by $\vee_2$, $\wedge_2$. Hence its image are in the set of $\{a,b\}$, and by the property of these operators, $f|_{\{a,b\}^n}$ is monotonic and reproducing, that is, (adopting the notations in (\ref{4.01})) $\forall \tilde{x}_1,\cdots,\tilde{x}_n\in\{0,1\}$,
\begin{align*}
	&f_{ab}(\tilde{x}_1,\cdots,\tilde{x}_n)\in\{a,b\},\\
	&f_{ab}(\tilde{x}_1,\cdots,\tilde{x}_{i-1},1,\tilde{x}_{i+1},\cdots,\tilde{x}_n)-f_{ab}(\tilde{x}_1,\cdots,0,\cdots,\tilde{x}_n)\in\{a-b,0\},\quad\forall i=1,\cdots,n,\\
	&f_{ab}(1,\cdots,1)=a,~f_{ab}(0,\cdots,0)=b,
\end{align*}
where in the second equation the variables only differ in the $i$-th entry, $i=1,\cdots,n$, and it is equivalent to say that $f(x_1,\cdots,x_n)\leqslant f(y_1,\cdots,y_n)$ if $x_i\leqslant y_i$, $x_i, y_i\in\{a,b\}$, $i=1,\cdots,n$. Translating the above equations into the vector expression yields 
\begin{align*}
	&M_{ab}\tilde{x}_1\cdots\tilde{x}_n\in\{\delta_k^{i_a},\delta_k^{i_b}\};\\
	&M_{ab}\tilde{x}_1\cdots\delta_2^1\cdots\tilde{x}_n-M_{ab}\tilde{x}_1\cdots\delta_2^2\cdots\tilde{x}_n=M_{ab}W_{[2,2^{i-1}]}(\delta_2^1-\delta_2^2)\tilde{x}_1\cdots\tilde{x}_{i-1}\tilde{x}_{i+1}\cdots\tilde{x}_n\in\{\delta_k^{i_a}-\delta_k^{i_b},\mathbf{0}\};\\
	&M_{ab}\delta_2^1\cdots\delta_2^1=\delta_k^{i_a}, ~M_{ab}\delta_2^2\cdots\delta_2^2=\delta_k^{i_b},\\
	&\forall \tilde{x}_1,\cdots,\tilde{x}_n\in\Delta_2,~i=1,\cdots,n,
\end{align*}
which are (\ref{eq:s1})(\ref{eq:s2})(\ref{eq:s3}) respectively.

(Sufficiency) Conversely, if $f|_{\{a,b\}^n}$ satisfies the conditions (\ref{eq:s1})(\ref{eq:s2})(\ref{eq:s3}), it can be viewed as generated by $\vee_2$, $\wedge_2$ over the lattice $\{a,b\}$, since according to the classical result by E. Post \cite{po41}, the class of reproducing and monotonic $n$-ary functions over ${\mathcal D}_2$ are exactly the one generated by $\vee_2$, $\wedge_2$. Therefore $a$, $b$ has the order relation obtained from these operators.
\hfill $\Box$

Theorem \ref{t4.2} shows that given any lattice function $f$ over a $k$-element lattice, we can find all the comparable pairs $a,b$ in ${\mathcal D}_k$ by checking the conditions on corresponding $M_{ab}$, and further obtain an undirected graph of $k$ nodes, called the comparability graph, where there is an edge between two nodes if and only if they are comparable. 

Then algorithms can be applied to the comparability graph, assigning to it an orientation to recover the partial order on ${\mathcal D}_k$, which is, however, not unique \cite{bra99}.  Details of these algorithms will be stated later in Section \ref{S4.2}.

\begin{remark}
	The significance of Theorem \ref{t4.2} is that, once the underlying lattice structure of a network over a lattice is recovered, one may use it to simplify the analysis of the control problems. For example, if we find out that the network is defined over a product lattice, then the results in Section \ref{S3B} can be applied to decompose the network.
\end{remark}

\subsection{Realization of the Underlying Lattice Structure of a Network}\label{S4.2}

Next we consider designing a lattice structure over ${\mathcal D}_k$ to make a $k$-valued network be generated by some basic operators of the lattice.  The only data available is the structure matrix. The construction requires two steps: first, derive a relation from $f$ according to Theorem \ref{t4.2}; second, verify if it is a partial order making ${\mathcal D}_k$ a lattice.

According to the results in \cite{lau06} (chapter 11, section 11.4), monotonic functions over a lattice are generated by $\vee$, $\wedge$ and piecewise constant functions $m_{a,b}$, defined as
\begin{align}\label{4.3}
	m_{a,b}(x)=\begin{cases}
		b,\quad a\leqslant x\\
		\mathbf{0}, \quad \text{otherwise}.
	\end{cases}
\end{align}
We call $\vee$, $\wedge$, $\{m_{a,b}\}_{a,b\in {\mathcal D}_k}$ basic operators on a lattice.

We present the following algorithm to construct a lattice structure on ${\mathcal D}_k$ to make the network generated by basic operators on this lattice. 

\begin{proposition}\label{p4.3}
    Consider a network over ${\mathcal D}_k$ with $n$ nodes and its ASSR as $x(t+1)=Mx(t)$, where $x(t)=x(t)=\ltimes_{j=1}^nx_j(t)$, $M\in{\mathcal L}_{k^n\times k^n}$. Execute the following procedure to derive a partial order on ${\mathcal D}_k$:
	\begin{itemize}
		\item Step 1: Apply Proposition \ref{p4.1} to decompose $M$ into component-wise form $M_1,\cdots,M_n$, where $M_i$ is the structure matrix of the $i$-th node $x_i(t+1)=f_i(x_1(t),\cdots,x_n(t))$ with ASSR $x_i(t+1)=M_ix(t)$, $i=1,\cdots,n$. {Remove the dumb indices in each function $M_i$ (see the algorithm provided in \cite{che10a}).}
		\item  Step 2: For the $i$-th function $f_i$ with structure matrix $M_i$, draw the corresponding graph $G_i$ of $n$ nodes, where there is an edge between nodes $a, b\in{\mathcal D}_k$ if and only if $(a,b)$ is a pair on which $f_i$ satisfies the conditions (\ref{eq:s1})(\ref{eq:s2}). Take the conjunction of the graphs $G_i$, $i=1,\cdots,n$, denote it by $G$.
		\item Step 3: Give an orientation to the resulting graph $G$, verify whether it is transitive. If it is, then the topological sorting of the oriented graph  derives a partial order on ${\mathcal D}_k$.
	\end{itemize}
 
	If the partial order derived from the network over ${\mathcal D}_k$ following the above procedure makes a lattice, then the network is generated from the class of functions $\{\vee,\wedge,m_{a,b}\}$, where $\vee$, $\wedge$ are the $\sup$, $\inf$ operators of the lattice respectively, and $m_{a,b}$ are piecewise constant functions defined as in (\ref{4.3}), $a,b\in{\mathcal D}_k$.
\end{proposition}

\noindent {\it Proof.} The pairs on which the function with structure matrix $M_i$ is monotonic (i.e. satisfies (\ref{eq:s1})(\ref{eq:s2}), by the proof of Theorem \ref{t4.2}) are candidates of comparable pairs to allow this function to be generated from basic operators of a lattice (see chapter 11 of \cite{lau06}). After checking the updating function of the $i$-th node on each pair $\{a,b\}\subset{\mathcal D}_k$, we obtain the comparability graph of the partial order for the dynamics of the $i$-th node. Since only when two elements are comparable in the order relation derived from each node can they be labelled as comparable in the partial order corresponding to the whole system, one need to take conjunction of the $n$ graphs (that is, there is an edge between $\{a,b\}$ in graph $G$ if and only if such edge exists for each graph $G_i$, $i=1,\cdots,n$), and check the conditions for being a lattice on the corresponding order relation after orienting the graph.
\hfill$\Box$

Proposition \ref{p4.3} gives a criterion for whether a network can be realized through basic operators of some lattice when only its structure matrix is known.

\begin{remark}
    The  procedure of constructing the graph $G$ in Step 2 is depicted in pseudo-code form as in the following table \ref{alg1}. Let $(N,E)$ be an undirected graph of $k$ nodes, i.e. $N:=\{1,\cdots,k\}$, and $E:=\{e(a,b):{\mathcal D}_k\times{\mathcal D}_k\rightarrow\{0,1\}\}_{a,b\in1,\cdots,k}$, where $e(a,b)=1$ if and only if there is an edge between the vertices $a$, $b$. Then Algorithm \ref{alg1} returns the adjacency information of the vertices.
\end{remark}

\begin{algorithm}[htbp]
\caption{Deriving the comparability graph from the ASSR}\label{alg1}
\hspace*{0.02in} {\bf Input:}
Structure matrices $M_i\in{\mathcal L}_{k\times k^n}$, $i=1,\cdots,n$.\\
\hspace*{0.02in} {\bf Output:}
Graph $G=(N,E)$.
\begin{algorithmic}[1]
\STATE{Set $e(a,b)=0$, $e_i(a,b)=0$ for all $a,b\in\{1,\cdots,k\}$, $i=1,\cdots,n$;}
\FOR{$i=1:n$}
\FOR{$a=1:k$}
\FOR{$b=1:a$}
\IF {$e_i(a,b)=0$}
\STATE {{\bf goto} here;}
\ELSE
\STATE{${M}^{ab}_i=M_i\prod_{j=0}^{n-1}(I_{2^j}\otimes \delta_k[i_a,i_b])$;}
\IF {$M_i^{ab}$ satisfies (\ref{eq:parent})}
\STATE{$e_i(a,b)=1$;}
\ELSE
\STATE{$e_i(a,b)=0$;}
\ENDIF
\ENDIF
\STATE {here}
\ENDFOR
\ENDFOR
\STATE{$e(a,b)=e(a,b)\wedge e_i(a,b)$;}
\ENDFOR
\STATE\RETURN $N=\{1,\cdots,k\}$, $E=\{e(a,b)\}_{a,b=1,\cdots,k}$
\end{algorithmic}
\end{algorithm}

\begin{remark}\label{r4.5}
	The algorithm in Step 3 for transitive orientation have been well established since 1999 \cite{bra99}, and we only give a sketch of it due to the limitation of space. Let $G$ be a prime undirected graph with vertex set $V$, start with a partition $\{\{v\},V\backslash\{v\}\}$ and define it as an ordered list. Then choose a pivot vertex $x\in V$, and split each partition class $Y$ into two parts: the vertices adjacent to $x$ (denoted by $Y_a$) and those non-adjacent to $x$ (denoted by $Y_n$) and place them in consecutive positions in the ordered list. Let $Y_n$ occupy the earlier in the two new positions if $x$ proceeds $Y$, and the latter if not. Then move on to another pivot vertex. Refining the partition through the procedure of choosing pivots, one will obtain an ordered partition where each class only contains a single element, which is the topological sorting of the orientation. One may refer to \cite{mc94,mc99} for proofs and details.
\end{remark}

In the end, we claim that there exist canonical expressions for monotonic functions over a lattice.

\begin{proposition}[\cite{lau06}, Section 11.4]\label{p4.5}
	If $f:{\mathcal D}_k^n\rightarrow{\mathcal D}_k$ is monotonic over a lattice $({\mathcal D}_k,\leqslant)$. Then $\forall \bm{x}=(x_1,\cdots,x_n)\in{\mathcal D}_k^n$,
	\begin{align}\label{4.5}
		f(\bm x)=\bigvee_{\bm{a}=(a_1,\cdots,a_n)\in{\mathcal D}_k^n}\Big(\bigwedge_{i=1}^nm_{a_i,f(\bm{a})}(x_i)\Big),
	\end{align}
	where $m_{a_i,f(\bm{a})}$ is defined as in (\ref{4.3}).
\end{proposition}

One can see that additional constraints are needed to make sure that a function is a lattice function (that is, generated solely by $\vee$, $\wedge$); however, finding such conditions would be difficult due to (\ref{4.5}), since any lattice function is monotonic and hence will always allow an expression that contains $m_{ab}$.

\section{Numerical Examples}\label{S5}

This section provides numerical examples to illustrate the results in sections \ref{S3} and \ref{S4}. 

We first consider the observability problem of a control network defined over a product lattice, and see how it can be solved through decomposing the system into subsystems over factor lattices.

\begin{example}\label{ex5.1}
	Consider a network over the lattice ${\mathcal D}_2\times{\mathcal D}_3$, where ${\mathcal D}_2=\{0,1\}$, ${\mathcal D}_3=\{0,1,2\}$ are ordered as chains canonically. The Hasse diagram of the product lattice is shown in Figure \ref{Figpa.4.0}.
	
	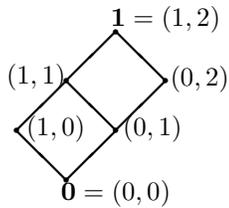
\begin{figure}
		\centering
		\setlength{\unitlength}{0.33 cm}
		\begin{picture}(8,8)\thicklines
			\put(1,3){\circle*{0.2}}
			\put(3,1){\circle*{0.2}}
			\put(3,5){\circle*{0.2}}
			\put(5,3){\circle*{0.2}}
			\put(5,7){\circle*{0.2}}
			\put(7,5){\circle*{0.2}}
			\put(1,3){\line(1,-1){2}}
			\put(1,3){\line(1,1){2}}
			\put(3,1){\line(1,1){2}}
			\put(3,5){\line(1,-1){2}}
			\put(3,5){\line(1,1){2}}
			\put(5,3){\line(1,1){2}}
			\put(5,7){\line(1,-1){2}}
			\put(2.8,0.2){\small${\bf 0}=(0,0)$}
			\put(1.4,2.8){\small$(1,0)$}
			\put(5.3,2.8){\small$(0,1)$}
			\put(0.6,4.9){\small$(1,1)$}
			\put(7.2,4.8){\small$(0,2)$}
			\put(4.8,7.2){\small${\mathbf 1}=(1,2)$}
		\end{picture}
		\caption{Hasse diagram of ${\mathcal D}_2\times{\mathcal D}_3$\label{Figpa.4.0}}
	\end{figure}
	
	The dynamics of the network is
	\begin{align}\label{4.0}
		\begin{cases}
			x_1(t+1)=x_2(t)\wedge u(t)\\
			x_2(t+1)=x_1(t)\vee x_2(t)\\
		\end{cases},
		\quad y(t)=x_2(t).
	\end{align}
	Converting (\ref{4.0}) to ASSR, we obtain the structure matrices of the dynamics and the output respectively as
	\begin{align*}
		L&=\delta_{36}[1,3,3,5,5,\cdots,32,33,34,35,36]\in{\mathcal L}_{36\times 216},\\
		E&=\delta_6[1,2,3,4,5,6,\cdots,1,2,3,4,5,6]\in{\mathcal L}_{6\times 36}.
	\end{align*}
	We will consider the observability of the system (\ref{4.0}), showing how existing results can be applied to the ASSR, and how it could be solved through its subnetworks over sublattices. Suppose $x_1(t)=(x^1_1(t),x^1_2(t))$, $x_2(t)=(x^2_1(t),x^2_2(t))$, $y(t)=(y_1(t),y_2(t))\in\Delta_2\times\Delta_3$, then one can calculate the ASSR of the corresponding subnetworks over the factor rings respectively as
	\begin{align*}
		\begin{array}{rl}
			\begin{cases}
				z_1(t+1)=\tilde{L}_1z_1(t),\\
				y_1(t)=E_1z_1(t);
			\end{cases}
			\begin{cases}
				z_2(t+1)=\tilde{L}_2z_2(t),\\
				y_2(t)=E_2z_2(t),
			\end{cases}
		\end{array}
	\end{align*}
	where $z_1(t)=x^1_1(t)x^1_2(t)$, $z_2(t)=x^1_2(t)x^2_2(t)$, and
	\begin{align*}
		\begin{array}{l}
			\tilde{L}=\delta_4[1,2,2,2,1,2,4,4],\quad
			E_1=\delta_2[1,2,1,2],\\
			\tilde{L}_2=\delta_9[1,2,3,2,2,3,3,3,3,1,2,3,5,5,6,6,6,6,1,2,3,5,5,6,9,9,9],\\
			E_2=\delta_3[1,2,3,1,2,3,1,2,3].
		\end{array}
	\end{align*}
	First consider the network over ${\mathcal D}_2$. Following \cite{che18}, construct an auxiliary network as
	\begin{align*}
		\begin{cases}
			z_1(t+1)=\tilde{L}_1u(t)z_1(t)\\
			z^*_1(t+1)=\tilde{L}_1u(t)z^*_1(t),
		\end{cases}
	\end{align*}
	set $w(t)=z(t)z^*(t)$, then the ASSR is
	$$
	w(t+1)=Gu(t)w(t),
	$$
	where $G=\tilde{L}_1(I_8\otimes\tilde{L}_1)(I_2\otimes W_{[2,4]})\mathrm{diag}(\delta_2^1,\delta_2^2)$.
	
	Solving the output distinguishable pairs, we derive the set
	$$
	W=\{\delta_{16}^2,\delta_{16}^4,\delta_{16}^5,\delta_{16}^7,\delta_{16}^{10},\delta_{16}^{12},\delta_{16}^{13},\delta_{16}^{15}\},
	$$
	Then the controllability matrix from $\Delta_{16}$ to $W$ is 
	\begin{align*}
		C_W^1&=I_W[\sum_{i=1}^{16}(G\delta_2^1+G\delta_2^2)^{(i)}]\\
		&=[0,1,1,1,1,0,0,0,1,0,0,0,1,0,0,0],
	\end{align*}
	where $I_W$ is the index matrix of $W$.
	
	Therefore the distinguishable pairs of subnetwork is $S_1=\{(\delta_4^1,\delta_4^2),(\delta_4^1,\delta_4^3),(\delta_4^1,\delta_4^4)\}$, the system is not observable.
	
	A similar argument can be applied to the subsystem over ${\mathcal D}_3$, showing that
	\begin{small}
		\begin{align*}
			\begin{array}{ccl}
				C_W^2&=&[0,1,1,1,1,1,1,1,1,1,0,1,0,0,1,1,1,1,1,1,0,1,1,0,0,0,0,1,0,1,0,0,1,1,1,1,1,0,1,0,\\
    &~&0,1,1,1,1,1,1,0,1,1,0,0,0,0,1,1,0,1,1,0,0,0,0,1,1,0,1,1,0,0,0,0,1,1,0,1,1,0,0,0,0],\\
			\end{array}
		\end{align*}
	\end{small}
	and the distinguishable pairs of the subnetwork is
	\begin{small}
		\begin{align*}
			\begin{array}{rcl}
				S_2&=&\{(\delta_9^1,\delta_9^2),(\delta_9^1,\delta_9^3),(\delta_9^1,\delta_9^4),
				(\delta_9^1,\delta_9^5),(\delta_9^1,\delta_9^6),(\delta_9^1,\delta_9^7),(\delta_9^1,\delta_9^8),(\delta_9^1,\delta_9^9),(\delta_9^2,\delta_9^3),(\delta_9^2,\delta_9^6),(\delta_9^2,\delta_9^7),(\delta_9^2,\delta_9^8),\\
				&~&(\delta_9^2,\delta_9^9),(\delta_9^3,\delta_9^4),(\delta_9^3,\delta_9^5),(\delta_9^4,\delta_9^6),(\delta_9^4,\delta_9^7),(\delta_9^4,\delta_9^8),(\delta_9^4,\delta_9^9),(\delta_9^5,\delta_9^6),(\delta_9^5,\delta_9^7),(\delta_9^5,\delta_9^8),(\delta_9^5,\delta_9^9)\}
			\end{array}
		\end{align*}
	\end{small}
	Therefore we conclude that the whole system is not observable. To be precise, $(\delta_2^{i_1},\delta_2^{j_1}), (\delta_2^{i_2},\delta_2^{j_2})\in ({\Delta}_2\times{\Delta}_3)^2$ is distinguishable if and only if, $(\delta_2^{i_1},\delta_2^{j_1})\in S_1$, $(\delta_2^{i_2},\delta_3^{j_2})\in S_2$.
	
\end{example}

Next, we give an example of how the algorithm in Proposition \ref{p4.3} works to construct a lattice structure on a finite set from the structure matrix of a network, so that this network can be expressed through basic operators on the lattice.

\begin{example}\label{ex5.2}
	Consider the following network over ${\mathcal D}_5$ with $3$ nodes and ASSR as 
	\begin{align}\label{5.4}
		x(t+1)=Mx(t),
	\end{align}
	where $x(t)=x_1(t)x_2(t)x_3(t)$, and
	\begin{small}
		\begin{align*}
			\begin{array}{ccl}
				M=&\delta_{125}& [1,31,61,91,121,26,31,61,91,121,51,56,61,91,121,76,81,86,91,121,101,106,111,116,121,\\
				~&        ~& 1,32,61,92,122,32,32,67,92,117,51,57,61,92,122,82,82,92,92,117,102,107,112,117,122,1,\\
				~&        ~& 31,63,93,123,26,31,63,93,123,63,68,63,93,118,88,93,88,93,118,103,108,113,118,123,1,\\
				~&        ~& 32,63,94,124,32,32,69,94,119,63,69,63,94,119,94,94,94,94,119,124,119,119,119,124,1,\\
				~&        ~& 32,63,94,125,27,32,64,94,125,53,59,63,94,125,99,94,94,94,125,125,120,120,120,125].
			\end{array}
		\end{align*}
	\end{small}
	We aim to check if there is an underlying partial order on ${\mathcal D}_5$ to make this system be one defined over a lattice, or to realize the transition law by some lattice operators.
	
	Decomposition of $M$ according to Proposition \ref{p4.1}, we derive the structure matrix of the dynamics of each node as
	$$
	M_i={\mathbf 1}_{5^{i-1}}\otimes I_{5}\otimes{\mathbf 1}_{5^{3-i}}\in{\mathcal L}_{5\times 5^3},~i=1,2,3.
	$$
	Next, we try to derive comparability graphs from them. Checking the conditions (\ref{eq:s1})(\ref{eq:s2}) on $M_1$, $M_2$ and $M_3$ for each pair of $a,b\in{\mathcal D}_5$, one can obtain their comparability graphs, as shown in Figure \ref{Figpa.4.1}, \ref{Figpa.4.2} and \ref{Figpa.4.3} respectively.
	\begin{figure}
		\centering
		\subfloat[\label{Figpa.4.1}]{
			\setlength{\unitlength}{0.5cm}
			\begin{picture}(5,4)\thicklines
				\put(2,0.5){\circle*{0.2}}
				\put(4,0.5){\circle*{0.2}}
				\put(1,2.5){\circle*{0.2}}
				\put(5,2.5){\circle*{0.2}}
				\put(3,3.5){\circle*{0.2}}
				\put(2,0.5){\line(-1,2){1}}
				\put(1,2.5){\line(2,1){2}}
				\put(3,3.5){\line(2,-1){2}}
				\put(5,2.5){\line(-1,-2){1}}
				\put(5,2.5){\line(-1,0){4}}
				\put(2,0.5){\line(3,2){3}}
				\put(1,2.5){\line(3,-2){3}}
				\put(2,0.5){\line(1,3){1}}
				\put(3,3.5){\line(1,-3){1}}
				\put(2,0.5){\line(1,0){2}}
				\put(1.3,0.5){$2$}
				\put(0.5,2.5){$1$}
				\put(5.2,2.5){$4$}
				\put(3,3.7){$0$}
				\put(4.3,0.5){$3$}
		\end{picture}}
		\subfloat[\label{Figpa.4.2}]{
			\setlength{\unitlength}{0.5cm}
			\begin{picture}(5,4)\thicklines
				\put(2,0.5){\circle*{0.2}}
				\put(4,0.5){\circle*{0.2}}
				\put(1,2.5){\circle*{0.2}}
				\put(5,2.5){\circle*{0.2}}
				\put(3,3.5){\circle*{0.2}}
				\put(2,0.5){\line(-1,2){1}}
				\put(1,2.5){\line(2,1){2}}
				\put(3,3.5){\line(2,-1){2}}
				\put(5,2.5){\line(-1,-2){1}}
				\put(5,2.5){\line(-1,0){4}}
				\put(2,0.5){\line(3,2){3}}
				\put(1,2.5){\line(3,-2){3}}
				\put(1.3,0.5){$2$}
				\put(0.5,2.5){$1$}
				\put(5.2,2.5){$4$}
				\put(3,3.7){$0$}
				\put(4.3,0.5){$3$}
		\end{picture}}
		\subfloat[\label{Figpa.4.3}]{
			\setlength{\unitlength}{0.5cm}
			\begin{picture}(5,4)\thicklines
				\put(2,0.5){\circle*{0.2}}
				\put(4,0.5){\circle*{0.2}}
				\put(1,2.5){\circle*{0.2}}
				\put(5,2.5){\circle*{0.2}}
				\put(3,3.5){\circle*{0.2}}
				\put(2,0.5){\line(-1,2){1}}
				\put(1,2.5){\line(2,1){2}}
				\put(3,3.5){\line(2,-1){2}}
				\put(5,2.5){\line(-1,-2){1}}
				\put(5,2.5){\line(-1,0){4}}
				\put(2,0.5){\line(3,2){3}}
				\put(1,2.5){\line(3,-2){3}}
				\put(2,0.5){\line(1,3){1}}
				\put(3,3.5){\line(1,-3){1}}
				\put(1.3,0.5){$2$}
				\put(0.5,2.5){$1$}
				\put(5.2,2.5){$4$}
				\put(3,3.7){$0$}
				\put(4.3,0.5){$3$}
		\end{picture}}
		\caption{Comparability graphs corresponding to $M_1$, $M_2$, $M_3$}
	\end{figure}
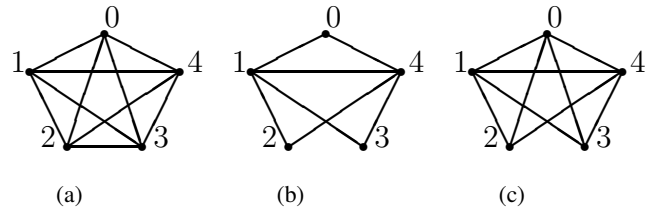
	Then take their conjunction to get the comparability graph of the system (which is, naturally, the one in Figure \ref{Figpa.4.2}) and orient it according to the vertex-partitioning algorithm in \cite{mc94}. First, choose $1$ as a source and the initial partition is $\{\{1\},\{2,3,4,0\}\}$. After a pivot on $2$, since $\{3,0\}$ is not adjacent to $2$, it goes before $\{2,4\}$; then after a pivot on $3$, we obtain the topological sorting of the graph as
	\begin{align}\label{5.2}
		1~|~0~|~3~|~2~|~4
	\end{align}
	where an element is less than the ones following it in the sequence. Applying (\ref{5.2}) to Figure \ref{Fig1} gives the conjunction of comparability graphs a transitive orientation, as shown in Figure \ref{Fig2}. Hence we obtain a partial order on ${\mathcal D}_5$, denoted by $R$.
	
	\begin{figure}
		\centering
		\subfloat[\label{Fig1}]{
			\setlength{\unitlength}{0.5cm}
			\begin{picture}(5,4)\thicklines
				\put(2,0.5){\circle*{0.2}}
				\put(4,0.5){\circle*{0.2}}
				\put(1,2.5){\circle*{0.2}}
				\put(5,2.5){\circle*{0.2}}
				\put(3,3.5){\circle*{0.2}}
				\put(1,2.5){\vector(1,-2){0.95}}
				\put(1,2.5){\vector(2,1){1.9}}
				\put(3,3.5){\vector(2,-1){1.9}}
				\put(4,0.5){\vector(1,2){0.9}}
				\put(1,2.5){\vector(1,0){3.8}}
				\put(2,0.5){\vector(3,2){2.9}}
				\put(1,2.5){\vector(3,-2){2.9}}
				\put(1.3,0.5){$2$}
				\put(0.5,2.5){$1$}
				\put(5.2,2.5){$4$}
				\put(3,3.7){$0$}
				\put(4.3,0.5){$3$}
		\end{picture}}
		\quad
		\subfloat[\label{Fig2}]{
			\setlength{\unitlength}{0.7 cm}
			\begin{picture}(3,2.2)\thicklines
				\put(2,0.5){\circle*{0.15}}
				\put(1,1.5){\circle*{0.15}}
				\put(3,1.5){\circle*{0.15}}
				\put(2,1.5){\circle*{0.15}}
				\put(2,2.5){\circle*{0.15}}
				\put(2,0.5){\line(-1,1){1}}
				\put(2,0.5){\line(0,1){2}}
				\put(2,0.5){\line(1,1){1}}
				\put(2,2.5){\line(-1,-1){1}}
				\put(2,2.5){\line(1,-1){1}}
				\put(2.3,0.3){$1$}
				\put(0.6,1.5){$2$}
				\put(2.2,1.5){$3$}
				\put(3.2,1.5){$0$}
				\put(2,2.6){$4$}
		\end{picture}}
		\caption{Transitive orientation of the conjunction of comparability graphs and the Hasse diagram of the lattice\label{Figpa.4.5}}
	\end{figure}
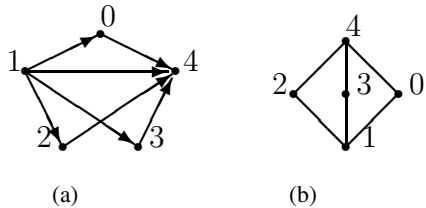
	
	Finally, after checking the conditions on this orientation for being a lattice over ${\mathcal D}_5$, we may draw the Hasse diagram of this order relation, shown as in Figure \ref{Figpa.4.5}.

	Therefore, $({\mathcal D}_5,R)$ is the lattice structure that makes the network (\ref{5.4}) become one generated by basic operators on the lattice, and one can write it explicitly in these operators according to (\ref{4.5}).

\end{example}

\section{Conclusions}\label{S6}

In this paper control networks over finite lattices are firstly considered. Using STP, the ASSR of such networks is obtained, making control problems easily solvable. Algorithms are proposed to recover the partial order from a lattice function, and construct lattice structures for an arbitrary network to allow it to be generated by basic operators over the lattice.

In applications, the results in Section \ref{S3B} and Section \ref{S4} can be combined: given a network over an unknown lattice, one may first recover the order relation using Theorem \ref{t4.2}, then decompose the system into subsystems defined over factor sublattices if the lattice is a product lattice, and Theorem \ref{t3} will help to simplify the computation of the control properties. {For example, if we apply the algorithm in Proposition \ref{p4.3} to the network dynamics in Example \ref{ex5.1}, we will reconstruct the oriented comparability graph as in Figure \ref{f5}, and further can recover the lattice structure depicted as in Figure \ref{Figpa.4.0}. Then the decomposition can be applied, reducing the complexity of the system analysis from $O(6^2)$ to $O(2^2+3^2)$.}
    \begin{figure}
		\centering
			\setlength{\unitlength}{0.5cm}
			\begin{picture}(6,6)\thicklines
				\put(0,2){\circle*{0.2}}
				\put(6,2){\circle*{0.2}}
				\put(0,4){\circle*{0.2}}
				\put(6,4){\circle*{0.2}}
				\put(3,6){\circle*{0.2}}
                \put(3,0){\circle*{0.2}}
				\put(3,0){\vector(-3,2){3}}
                \put(3,0){\vector(3,2){2.9}}
				\put(3,0){\vector(-3,4){3}}
                \put(3,0){\vector(3,4){3}}
				\put(3,0){\vector(0,1){5.9}}
				\put(0,2){\vector(1,0){5.65}}
				\put(0,2){\vector(3,1){6}}
                \put(0,2){\vector(3,4){2.7}}
				\put(0,4){\vector(3,-1){5.9}}
				\put(0,4){\vector(3,2){3}}
                \put(3,6){\vector(3,-4){2.85}}
                \put(6,4){\vector(0,-1){1.9}}
				\put(-0.3,1.3){$2$}
				\put(3.2,-0.5){$1$}
				\put(-0.3,4.2){$3$}
				\put(3.2,6){$4$}
				\put(5.8,4.2){$5$}
                \put(5.9,1.3){$6$}
            \end{picture}
            \caption{Oriented comparability graph over ${\mathcal D}_6$ in Example \ref{ex5.1}}\label{f5}
    \end{figure}

Although due to the canonical expression (\ref{4.5}), the difficulty of finding restrictions on the monotonic function class over a lattice to exclude the functions $m_{a,b}$ seems intrinsic, a criteria for a multiple-valued logical function to be generated solely by the operators $\vee$ and $\wedge $ is still a problem worth investigating.

Meanwhile, using the framework proposed in this paper, problems concerning switched, delayed and probabilistic networks over finite lattices can be further investigated.

\Acknowledgements{This work is supported partly by the National Natural Science Foundation of China (NSFC) under Grants 62073315 and 62350037.}

\end{document}